\documentclass[11pt,letterpaper,onecolumn]{article}

\usepackage{amsmath}
\usepackage{amssymb}
\usepackage{exscale}
\usepackage{amsthm}
\usepackage{epsfig}
\usepackage{setspace}

\theoremstyle{plain}

\newtheorem{theorem}{Theorem}

\newtheorem{lemma}[theorem]{Lemma}

\newcommand{\setmid}{\,|\,}
\newcommand{\N}{{\mathbb{N}}}

\newcommand{\Z}{{\mathbb{Z}}}
\newcommand{\R}{{\mathbb{R}}}

\newcommand{\ee}{\varepsilon}

\newcommand{\la}{\left|}
\newcommand{\ra}{\right|}

\newcommand{\eps}{\epsilon}

\newcommand{\ben}{\begin{enumerate}}
\newcommand{\een}{\end{enumerate}}
\newcommand{\bec}{\begin{center}}
\newcommand{\ece}{\end{center}}

\newcommand{\beq}{\begin{equation}}
\newcommand{\eeq}{\end{equation}}
\newcommand{\ARR}{\textsc{arr}}

\newcommand{\ODD}{\textsc{odd}}
\newcommand{\INF}{\textsc{inf}}
\newcommand{\CON}{\textsc{con}}
\newcommand{\tmax}{t_{\max{}}}
\newcommand{\tpmax}{t^\prime_{\max{}}}
\newcommand{\DISC}{\textsc{disc}}

\begin{document}
\thispagestyle{empty}
\author{Joshua Cooper\thanks{ETH-Z\"urich, Institute of Theoretical Computer Science, \tt{jcooper@inf.ethz.ch}}
  \and Benjamin Doerr\thanks{Max--Planck--Institut f\"ur Informatik,
  Saarbr\"ucken, Germany}
  \and Joel Spencer\thanks{Courant Institute of Mathematical Sciences,
  New York, U.S.A., \tt{spencer@cims.nyu.edu}}
  \and G\'abor Tardos\thanks{Simon Fraser University, Canada and R\'enyi Institute, Hungary, \tt{tardos@cs.sfu.ca}}
}
\date{}

\title{Deterministic Random Walks on the Integers\thanks{The authors
    enjoyed the hospitality, generosity and the strong coffee of the
    R\'enyi Institute (Budapest) while doing this research.  Spencer's
    research was partially supported by EU Project Finite Structures
    003006; Doerr's by EU Research Training Network COMBSTRU; Cooper's
    by an NSF Postdoctoral Fellowship (USA, NSF Grant DMS-0303272);
    and Tardos's by the Hungarian National Scientific Research Fund
    grants OTKA T-046234, AT-048826 and NK-62321.}}

\maketitle

\begin{abstract}
  Jim Propp's $P$-machine, also known as the `rotor router model' is a
  simple deterministic process that simulates a random walk on a
  graph. Instead of distributing chips to randomly chosen neighbors,
  it serves the neighbors in a fixed order. 
  
  We investigate how well this process simulates a random walk. For
  the graph being the infinite path, we show that, independent of the
  starting configuration, at each time and on each vertex, the number
  of chips on this vertex deviates from the expected number of chips
  in the random walk model by at most a constant $c_1$, which is
  approximately $2.29$. For intervals of length $L$, this improves to
  a difference of $O(\log L)$, for the $L_2$ average of a contiguous
  set of intervals even to $O(\sqrt{\log L})$. All these bounds are tight. 
\end{abstract}

\section{The Propp Machine}

The following deterministic process was suggested by Jim Propp as an
attempt to derandomize random walks on infinite grids $\Z^d$: 

{\bf Rules of the Propp machine:} Each vertex $x \in \Z^d$ is
associated with a `rotor' and a cyclic permutation of the $2d$
cardinal directions of $\Z^d$. Each vertex may hold an arbitrary number
of `chips'. In each time step, each vertex sends out all its chips to
neighboring vertices in the following manner: The first chip is sent
into the direction the rotor is pointing, then the rotor direction is
updated to the next direction in the cyclic ordering. The second chip
is sent in this direction, the rotor is updated, and so on. As a result, the chips are distributed highly evenly among the
neighbors.

This process has attracted considerable attention recently. It turns out that
the Propp machine in several respects is a very good simulation of a
random walk. Used to simulate  internal diffusion limited aggregation
(repeatedly, a single chip is inserted at the origin, performs a rotor router walk
until it reaches an unoccupied position and occupies it), it was shown
by Levine and Peres~\cite{levineperes} that this derandomization
produces results that are extremely close to what a random walk would
have produced. See also Kleber's paper~\cite{kleber}, which adds
interesting experimental results: Having inserted three million chips,
the closest unoccupied site is at distance 976.45, the farthest
occupied site is at distance 978.06. Hence the occupied sites almost
form a perfect circle!

In \cite{cooperpropp,cooperproppx}, the authors consider the following
question: Start with an arbitrary initial position (that is, chips on
vertices and rotor directions), run the Propp machine for some time
and compare the number of chips on a vertex with the expected number
of chips a random walk run for the same amount of time would have
placed on that vertex. Apart from a technicality, which we defer to
the end of Section~\ref{secresults}, the answer is astonishing: For
any grid $\Z^d$, this difference (discrepancy) can be bounded by a constant,
independent of the number of chips, the run-time, the initial rotor
position and the cyclic permutation of the cardinal directions.

In this paper, we continue this work. We mainly regard the
one-dimensional case, but as will be visible from the proofs, our
methods can be extended to higher dimensions as well. Besides making
the constant precise (approximately 2.29), we show that the
differences become even better for larger intervals (both in space and
time). We also present a fairly general method to prove lower bounds
(the `arrow forcing theorem'). This shows that all our upper bounds
are actually sharp, including the aforementioned constant.

Instead of talking about the expected number of chips the random walk
produces on a vertex, we find it more convenient to think of the
following `linear' machine. Here, in each time step each vertex sends out
exactly the same (possibly non-integral) number of chips to each
neighbor. Hence, for a given starting configuration, after $t$
time-steps the number of chips in the linear model is exactly the
expected number of chips in the random walk model.

\section{Our Results}\label{secresults}

We obtain the following results (again, see the end of the section for
a slight technical restriction): Fix any starting configuration, that
is, the number of chips on each vertex and the position of the rotor
on each vertex. Now run both the Propp machine and the linear machine
for a fixed number of time-steps. Looking at the resulting chip
configurations, we have the following:
\begin{itemize}
\item On each vertex, the number of chips in both models deviates by
  at most a constant  $c_1 \approx 2.29$. One may interpret this to mean
  that the Propp machine simulates a random walk extremely well. In
  some sense, it is even better than the random walk. Recall that in a
  random walk a vertex holding $n$ chips only in expectation sends
  $n/2$ chips to the left and the right. With high probability, the
  actual numbers deviate from this by $\Omega(n^{1/2})$.
\item In each interval of length $L$, the number of chips that are in
  this interval in the Propp model deviates from that in the linear
  model by only $O(\log L)$ (instead of, e.g.,  $2.29 L$).
\item If we average this over all length $L$ intervals in some larger
  interval of $\Z$, things become even better. The average squared
  discrepancy in the length $L$ intervals also is only $O({\log L})$.
\end{itemize}

We may as well average over time. In the setting just fixed, denote by
$f(x,T)$ the sum of the numbers of chips on vertex $x$ in the last $T$
time steps in the Propp model, and by $E(x,T)$ the corresponding
number for the linear model. Then we have the following discrepancy
bounds:
\begin{itemize}
\item The discrepancy on a single vertex over a time interval of length
  $T$ is at most $|f(x,T) - E(x,T)| = O(T^{1/2})$. Hence a vertex
  cannot have too few or too many chips for a long time (it may,
  however, alternate having too few and too many chips and thus have an
  average $\Omega(1)$ discrepancy over time).
\item We may extend this to discrepancies in intervals in space and
  time: Let $I$ be some interval in $\Z$ having length $L$. Then the
  discrepancy in $I$ over a time interval of length $T$ satisfies
  \[\Big|\sum_{x \in I} f(x,T) - \sum_{x \in I} E(x,T)\Big| =
  \left\{\begin{array}{ll} O(L T^{1/2}) & \mbox{if } L
      \le2T^{1/2}, \\
      O(T \log(L T^{-1/2})) & \mbox{otherwise. } \end{array} \right.\]
  Hence if $L$ is small compared to $T^{1/2}$, we get $L$ times the
  single vertex discrepancy in a time interval of length $T$ (no
  significant cancellation in space); if $L$ is of larger order than
  $T^{1/2}$, we get $T$ times the $O(\log L)$ bound for intervals of
  length $L$ (no cancellation in time, the discrepancy cannot leave
  the large interval in short time).
\end{itemize}

All bounds stated above are sharp, that is, for each  bound there is a
starting configuration such that after suitable run-time of the
machines we find the claimed discrepancy on a suitable vertex, in
a suitable interval, etc.

\textbf{A technicality:} There is one limitation, which we only
briefly mentioned, but without which our results are not valid.
Note that since $\Z^d$ is a bipartite graph, the chips that
start on even vertices never mix with those which start on odd
positions. It looks as if we would play two games in one. This is not
true, however. The even chips and the odd ones may interfere with each
other through the rotors. Even worse, we may use the odd chips to
reset the arrows and thus mess up the even chips. Note that the odd
chips are not visible if we look at an even position after an even
run-time. An extension of the arrow-forcing theorem presented below
shows that we can indeed use the odd chips to arbitrarily reset the
rotors. This is equivalent to running the Propp machine in an
adversarial setting, where an adversary may decide each time where the
extra odd chips on a position is sent to. It is clear that in this
setting, the results above cannot be expected. We therefore assume
that \emph{the starting configuration has chips only on even
  positions} (``even starting configuration'') or only on odd positions (``odd
starting configuration''). An alternative, in
fact equivalent, solution would be to have two rotors on each vertex,
one for even and one for odd time steps. 

\section{The Basic Method}
\label{secbasic} 

For numbers $a$ and $b$ set $[a..b] = \{z \in \Z \setmid a \le z \le
b\}$ and $[b] = [1..b]$. For integers $m$ and $n$, we write $m \sim n$
if $m$ and $n$ have the same parity, that is, if $m - n$ is even.

For a fixed starting configuration, we use $f(x,t)$ to denote the
number of chips at time $t$ at position $x$ and $\ARR(x,t)$ to denote
the value of the arrow at time~$t$ and position $x$, i.e., $+1$ if it
points to the right, and $-1$ if it points to the left. We have:
\begin{eqnarray*}f(x,t+1)&=&f(x-1,t)/2+f(x+1,t)/2\\
&&+\ARR(x-1,t)(f(x-1,t)\bmod2)/2\\&&-\ARR(x+1,t)(f(x+1,t)\bmod2)/2,\\
\ARR(x,t+1)&=&(-1)^{f(x,t)}\ARR(x,t).\end{eqnarray*}

Note that after an even starting configuration if $x\sim t$ does not
hold, then we have $f(x,t)=0$ and $\ARR(x,t+1)=\ARR(x,t)$.

We consider the machine to be started at time $t=0$. Being a
deterministic process, the initial configuration (i.e., the values
$f(x,0)$ and $\ARR(x,0)$, $x \in \Z$) determines the configuration at
any time $t>0$ (i.e., the values $f(x,t)$ and $\ARR(x,t)$, $x \in
\Z$).  The totality of all configurations for $t > 0$ we term a {\it
  game}.  We call a configuration {\em even} if no chip is at an odd
position. Similarly, a position is {\em odd} if no chip is at an even
position. Clearly, an even position is always followed by an odd
position and vice versa.

By $E(x,t)$ we denote the expected number of chips on a vertex $x$
after running a random walk for $t$ steps (from the implicitly given
starting configuration). As described earlier, this is equal to the
number of chips on $x$ after running the linear machine for $t$
time-steps.

In the proofs, we need the following mixed notation.  Let
$E(x,t_1,t_2)$ be the expected number of chips at location $x$ and
time $t_2$ if a simple random walk were performed beginning from the
Propp machine's configuration at time $t_1$. In other words, this is
the number of chips on vertex $x$ after $t_1$ Propp and $t_2 - t_1$
linear steps.

Let $H(x,t)$ denote the probability that a chip arrives at location
$x$ at time $t\ge0$ in a simple random walk begun from the origin, i.e.,
$H(x,t) = 2^{-t} \binom{t}{(t+x)/2}$, if $t \sim x$, and $H(x,t) = 0$
otherwise. For $t>0$ let $\INF(y,t)$ denote the ``influence'' of a Propp step
of a single chip at distance $y$ with $t$ linear steps remaining
(compared to a linear step). More precisely, we compare the two 
probabilities that a chip on position $y$ reaches $0$ if (a) it is first
sent to the right (by a single Propp step) and then does a random walk
for the remaining $t-1$ time steps, or (b) it just does $t$ random
walk steps starting from $y$. Hence,
\[\INF(y,t) := H(y+1,t-1) - H(y,t).\] A simple calculation yields
\begin{equation}
  \label{eq:infcalc}
  \INF(y,t)=  - \tfrac{y}{t} H(y,t).
\end{equation}
This shows in particular, that $\INF(y,t)\le0$ for $y\ge0$ and $\INF(y,t)\ge0$
for $y\le0$. We have $\INF(0,t)=0$.

For notational convenience we extend the definitions of $H(x,t)$ and $\INF(x,t)$
by letting $H(x,t)=0$ for $t<0$ and $\INF(x,t)=0$ for $t\le0$.

Note that 
\begin{equation}
  \label{eq:inf}
  \INF(y,t) = \tfrac 12 H(y+1,t-1) - \tfrac 12 H(y-1,t-1) .
\end{equation}
Therefore, the first Propp step with arrow
pointing to the left has an influence of $-\INF(y,t)$. 

Using this notation, we can conveniently express the (signed)
discrepancy $f(x,t) - E(x,t)$ on a vertex $x$ using information about
when ``odd splits'' occurred.  It suffices to prove the result for the
vertex $x=0$.  Clearly, $E(0,t,t) = f(0,t)$ and $E(0,0,t) = E(0,t)$, so
that
  \begin{equation} \label{eq1}
    f(0,t) - E(0,t) = \sum_{s=0}^{t-1} \left ( E(0,s+1,t) - E(0,s,t) \right ).
  \end{equation}
  In comparing $E(0,s+1,t)$ and $E(0,s,t)$, note that whenever there
  are two chips on some vertex at time $s$, then these chips can be
  assumed to behave identically no matter whether the next step is a
  linear or a Propp step. Denote by $\ODD_{s}$ the set of locations
  which are occupied by an odd number of chips at time $s$. Then 
  \begin{eqnarray*}
    \lefteqn{E(0,s+1,t) - E(0,s,t)}\\ 
    & =&  \sum_{y \in \ODD_{s}} (H(y+\ARR(y,s),t-s-1) - H(y,t-s))\\ 
    & =& \sum_{y \in \ODD_{s}} \ARR(y,s) \, \INF(y,t-s).
  \end{eqnarray*}
  Therefore, appealing to (\ref{eq1}),
  $$
  f(0,t) - E(0,t) = \sum_{s=0}^{t-1} \sum_{y \in \ODD_{s}} \ARR(y,s) \, \INF(y,t-s).
  $$
Using $\INF(y,u)=0$ for $u\le0$ we can extend the summation above for all
non-negative integers $s$.

  Let $s_i(y)$ be the $i^\textrm{th}$ time that
  $y$ is occupied by an odd number of chips, beginning with $i = 0$.
Switching the order of summation and noting that the arrows flip
  each time there is an odd number of chips on a vertex, we have
  \begin{align}\label{eqn:disc}
    f(0,t) - E(0,t) & = \sum_{y \in \Z} \sum_{i \geq 0}
    \ARR(y,s_i(y)) \, \INF(y,t-s_i(y)) \nonumber \\
    & = \sum_{y \in \Z} \ARR(y,0) \sum_{i \geq 0} (-1)^i \, \INF(y,t-s_i(y)).
  \end{align}

This equation will be crucial in the remainder of the paper. It shows that
the discrepancy on a vertex only depends on the initial arrow
positions and the set of location-time pairs holding an odd number of
chips. 

In the remainder, we show that we can construct starting configurations
with arbitrary initial arrow positions and odd number of chips at
arbitrary sets of location-time pairs. This will be the heart of our 
lower bound proofs in the following sections. Here $\N_0$ denotes the set of
non-negative integers.

\begin{theorem}[Parity-forcing Theorem]
  For any initial position of the arrows and any
  $\pi: \Z \times \N_0 \to \{0,1\}$, there is an initial even
  configuration of the chips such that for all $x \in \Z$, $t \in
  \N_0$ such that $x \sim t$, $f(x,t)$ and 
  $\pi(x,t)$ have identical parity.
\end{theorem}

Since rotors change their direction if and only if the vertex has an
odd number of chips, the parity-forcing theorem is a consequence of the
following arrow-forcing statement.

\begin{theorem}[Arrow-forcing Theorem] 
  Let $\rho(x,t)\in \{-1,+1\}$ be arbitrarily defined for $t\ge0$
  integer and $x \sim t$.  Then there exists an even initial
  configuration that results in a game with $\ARR(x,t)=\rho(x,t)$ for
  all such $x$ and $t$. Similarly, if $\rho(x,t)$ is defined for $x\sim t+1$ a
  suitable odd initial configuration can be found.
\end{theorem}

\begin{proof}
By symmetry, it is enough to prove the first statement.

Assume the functions $f$ and $\ARR$ describe the game following an even initial
configuration, and for some $T\ge0$, we have $\ARR(x,t)=\rho(x,t)$ for all
$0\le t\le T+1$ and $x \sim t$. We modify the initial position by defining
$f'(x,0)=f(x,0)+\eps_x2^T$ for even $x$, while we have $f'(x,0)=0$ for odd
$x$ and $\ARR'(x,0)=\ARR(x,0)$ for all $x$. Here, $\eps_x\in \{0,1\}$ are to be
determined. 

Observe that a
pile of $2^T$ chips will split evenly $T$ times so that the
arrows at time $t\leq T$ remain the same. Our goal is to choose the values $\eps_x$
so that $\ARR'(x,t)=\rho(x,t)$ for $0\le t\le T+2$ and $x \sim t$.
As stated above this holds automatically for $t\le T$ as
$\ARR'(x,t)=\ARR(x,t)=\rho(x,t)$ in this case. For $t=T+1$ and $x-T-1$ even we
have $\ARR'(x,T+1)=\ARR'(x,T)=\ARR(x,T)=\ARR(x,T+1)=\rho(x,T+1)$ since we start
with an even configuration. To make sure the equality also
holds for $t=T+2$ we need to ensure that the parities of the piles $f'(x,T)$
are right. Observe that $\ARR'(x,T+2)=\ARR'(x,T)$ if $f'(x,T)$ is even,
otherwise $\ARR'(x,T+2)=-\ARR'(x,T)$. So for $x-T$ even we must make
$f'(x,T)$ even if and only if $\rho(x,T+2)=\rho(x,T)$. At time $T$ the
``extra" groups of $2^T$ chips have spread as in Pascal's Triangle and
we have
\[  f'(x,T) = f(x,T) + \sum_y \eps_y {T\choose {\frac{T+x-y}2}}
\]
where $x \sim T$ and the sum is over the even values of $y$ with
$|y-x|\le T$. As $f(x,T)$ are already given it suffices to set
the parity of the sum arbitrarily. For $T=0$ the sum is $\eps_x$ so this is
possible. For $T>0$ we express
\[ \sum_y \eps_y {T\choose {\frac{T+x-y}2}}=
\eps_{x+T} + h + \eps_{x-T}\]
where $h$ depends only on $\eps_{y}$ with $x-T<y< x+T$.
We now determine the $\eps_y$ sequentially.  We initialize by setting
$\eps_y=0$ for $-T<y\le T$.  The values $\eps_y$ for $y>T$
are set in increasing order. The value of $\eps_y$ is set so that the sum at
$x=y-T$ (and thus $f'(y-T,T)$) will have the correct parity. Similarly,
the values $\eps_y$ for $y\le-T$ are set in decreasing order. The value of
$\eps_y$ is set so that the sum at $x=y+T$ (and thus the $f'(y+T,T)$) will
have the correct parity.

Note that the above procedure changes an even initial configuration that
matches the prescription in $\rho$ for times $0\le t\le T+1$ into another even
initial configuration that matches the prescription in $\rho$ for times $0\le
t\le T+2$. We start by defining $f(x,0)=0$ for all $x$ (no chips anywhere) and
$\ARR(x,0)=\rho(x,0)$ for even $x$, while $\ARR(x,0)=\rho(x,1)$ for odd
$x$. We now have $\ARR(x,t)=\rho(x,t)$ for $0\le t\le1$ and $x \sim t$. We can
apply the above procedure repeatedly to get an even
initial configuration that satisfies the prescription in $\rho$ for an
ever increasing (but always finite) time period $0\le t<T$. Notice however,
that in the procedure we do not change the initial configuration of arrows
$\ARR(x,0)$ at all, and we change the initial number of chips $f(x,0)$ at
position $x$ only if $|x|\ge T$. Thus at any given position $x$ the initial
number of chips will be constant after the first $|x|$ iterations. This means
that the process converges to an (even) initial configuration. It is simple to
check that this limit configuration satisfies the statement of the theorem.
\end{proof}

\section{Discrepancy on a Single Vertex}\label{secvertex}

\begin{theorem} \label{thm1:z1constant} There exists a constant $c_1
  \approx 2.29$, independent of the initial (even) configuration, the
  time $t$, or the location $x$, so that
$$
|f(x,t) - E(x,t)| \leq c_1.
$$
\end{theorem}

The proof needs the following elementary fact. Let $X \subseteq
\R$. We call a
mapping $f: X \to \R$ \emph{unimodal}, if there is an $m \in
X$ such that $f$ is monotonically increasing in $\{x \in X \setmid x
\le m\}$ and $f$ is monotonically decreasing in $\{x \in X \setmid x \ge
m\}$. 

\begin{lemma}\label{lunimod1}
  Let $f: X \to \mathbb R$ be non-negative and unimodal.  Let $t_1, \ldots,
  t_n \in X$ such that $t_1 < \ldots < t_n$. Then
    \[\la\sum_{i = 1}^n (-1)^i f(t_i)\ra \le  \max_{x \in X} f(x).\]
\end{lemma}
\begin{proof}[Proof of Theorem~\ref{thm1:z1constant}]
  It suffices to prove the result for $x=0$. In case $t$ is even we start with
  an even configuration, if $t$ is odd, then with an odd configuration
  (otherwise both $f(0,t)$ and $E(0,t)$ would be zero with no discrepancy).
  
First we show that $\INF(y,u)$ with a fixed $y<0$ is a non-negative unimodal
function of $u$ if restricted to the values $u\sim y$. We have already
seen that it is non-negative. For the unimodality let $y<0$ and
$u>2$, $u\sim y$. We have
\begin{align*}
\INF(y,u) - \INF(y,u+2) & = - \frac{y}{u} H(y,u) + \frac{y}{u+2} H(y,u+2) \\
& = \frac{4+3u-y^2}{(u+2-y)(u+2+y)}\INF(y,u),
\end{align*}
whenever $u \geq y$.  Hence the difference is non-negative if $u\ge(y^2-4)/3$ and it is non-positive
if $u\le(y^2-4)/3$. Thus we have unimodality, with $\INF(y,u)$ taking its
maximum at the smallest value of $u$ exceeding $(y^2-4)/3$ with $u\sim
y$. Let $\tmax(y):=\lfloor(y^2-4)/3\rfloor$+2.  It is easy to check that $\tmax(y)\sim y$ always
holds, so we have that $\INF(y,u)$ takes its maximum for fixed $y<0$ at
$u=\tmax(y)$. For $y>0$ the values $\INF(y,u)$ are non-positive and by symmetry
the minimum is taken at $u=\tmax(y)$. For $y=0$ we have $\INF(y,u)=0$ for all
$u$. We have just proved the following:

\begin{lemma}\label{linfmax}
  For $y \in \mathbb{Z}$, the function $|\INF(y,t)|$ is maximized over
  all integers $t$ at $\tmax(y) = \lfloor (y^2-4)/3 \rfloor+2$.
\end{lemma}

To bound $|f(0,t)-E(0,t)|$ we use the formula (\ref{eqn:disc}) where the inner
sums are alternating sums, for which we can apply Lemma~\ref{lunimod1}, as
$y\sim t-s_i(y)$ holds by our even or odd starting position assumption. We get
\begin{align}
\nonumber |f(0,t) - E(0,t)| & \leq \sum_{y \in \Z} \bigg| \sum_{i \geq 0} (-1)^i \, \INF(y,t-s_i(y)) \bigg| \\
\nonumber & \leq \sum_{y \in \Z}  \max_u |\INF(y,u) |\\
\label{eq2} &= 2\sum_{y=1}^\infty|\INF(y,\tmax(y))|.
\end{align}

Here
\begin{align*}
|\INF(y,\tmax(y))|&=\frac
 y{\tmax(y)}2^{-\tmax(y)}\binom{\tmax(y)}{(\tmax(y)+y)/2}\\
&=O(y/(\tmax(y))^{3/2})=O(y^{-2}).
\end{align*}

and, therefore, (\ref{eq2}) implies that $|f(0,t)-E(0,t)|$ is bounded by
$$
c_1:=2\sum_{y=1}^\infty| \INF(y,\tpmax(y))| \approx 2.29,
$$
proving Theorem \ref{thm1:z1constant}.
\end{proof}

Amazingly, the constant $c_1$ defined above is best possible. Indeed, let $y>0$
be arbitrary and even and let $t_0=\tmax(y)$. We apply the Arrow-forcing Theorem
to find an even starting position that makes $\ARR(x,t)=-1$ if $x>0$ and
$t\le t_0-\tmax(x)$ or $x<0$ and $t>t_0-\tmax(x)$ and makes $\ARR(x,t)=-1$
otherwise. It is easy to verify that in this case at a position $|x|\le y$,
$x\ne0$ we have an odd number of chips exactly once at time $t_0-\tmax(x)$
and the formula (\ref{eqn:disc}) gives
$$f(0,t_0)-E(0,t_0)=2\sum_{x=1}^y|\INF(x,\tmax(x))|.$$

\section{Intervals in Space}\label{secsint}

In this section, we regard the discrepancy in intervals in $\Z$. For
an arbitrary finite subset $X$ of $\Z$ set
\begin{eqnarray*}
  f(X,t) &:=& \sum_{x \in X} f(x,t),\\
  E(X,t) &:=& \sum_{x \in X} E(x,t).
\end{eqnarray*}

We show that the discrepancy in an interval of length $L$ is $O(\log
L)$, and this is sharp. We need the following facts about $H$.

\begin{lemma}\label{lhuni}
  For all $x \in \Z$, $H(x,\cdot) : \{t \in \N_0 \setmid x \sim t\}
  \to \R; t \mapsto H(x,t)$ is unimodal. $H(x,t)$ is maximal for
  $t = x^2$. We have $H(x,x^2) = \Theta(|x|^{-1})$.
\end{lemma}

\begin{proof}
  Since $H(x,t-2) - H(x,t) = \frac{t-x^2}{t^2-t}
  H(x,t)$, we conclude that $H(x,t)$ is unimodal and for $|x|\ge2$
  it has exactly two maxima, namely $t = x^2 - 2$ and $t = x^2$, while for
  $|x|\le1$ the latter is the only maximum. A standard estimate gives the
  claimed order of magnitude. 
\end{proof}

\begin{theorem}\label{tspace}\label{intspace}
  For any even initial configuration, any time $t$ and any interval
  $X$ of length~$L$, \[|f(X,t) - E(X,t)| = O(\log L).\] For every $L>0$
  there is an
  even initial configuration, a time $t$ and an interval $X$ of length
  $L$ such that \[|f(X,t) - E(X,t)| = \Omega(\log L).\]
\end{theorem}

\begin{proof}
Using that the discrepancy of a single position is bounded we can assume
$X$ ends at an even position, and then by symmetry we may assume it ends
at $0$, i.e., $X = [-L+1..0]$. Fix any even
  initial configuration. By (\ref{eqn:disc}), we have \[f(X,t) -
  E(X,t) = \sum_{y \in \Z} \ARR(y,0) \sum_{x \in X} \sum_{i \ge 0}
  (-1)^i \INF(y-x,t-s_i(y)).\]
Note that the summation here can be restricted to values $x\sim t$, the
other values contribute zero.

  Let us call \[\CON(y) := \ARR(y,0) \sum_{x \in X} \sum_{i \ge 0}
  (-1)^i \INF(y-x,t-s_i(y))\] the contribution of the vertex $y$ to
  the discrepancy in the interval $X$. The
  contribution of a vertex depends on its distance from the
  interval~$X$. If $y$ is $\Omega(L)$ away from $X$, its influences on
  the various vertices of $X$ are roughly equal, and all such
  influences are quite small. In this case we bound its influence by
  $L$ times the one we computed in Theorem~\ref{thm1:z1constant}:
  
  Let $y>L$. By Lemmas~\ref{lunimod1} and~\ref{linfmax},
  \begin{align*}
    |\CON(y)| &= \bigg|\sum_{x \in X} \sum_{i \ge 0} (-1)^i
    \INF(y-x,t-s_i(y))\bigg|\\
    &\le \sum_{x \in X} \bigg|\sum_{i \ge 0} (-1)^i
    \INF(y-x,t-s_i(y))\bigg|\\
    &\le \sum_{x \in X} \max_t |\INF(y-x,t)|\\
    &\le O\bigg(\sum_{x \in X}  (y-x)^{-2}\bigg) = O(L y^{-2}).
  \end{align*}
  Hence the total contribution of these vertices is at most
  \begin{equation*}
    \sum_{y > L} |\CON(y)| = O\bigg(\sum_{y > L} L y^{-2}\bigg) = O(1)
  \end{equation*}
and by symmetry the same bound applies to the contribution of vertices
$y\le-2L$.

  We now turn to vertices $-2L<y\le L$. Here mainly those vertices of
  $X$ that are close to $y$ contribute to $\CON(y)$.  Hence, the
  approach above is too coarse. We use instead that (\ref{eq:inf}) yields a
  collapsing sum. To simplify our formulas we introduce
$$H'(x,t)=H(x-1,t)+H(x,t).$$
Note that $H'(x,t)=H(x,t)$ for $x\sim t$ and $H'(x,t)=H(x-1,t)$ otherwise.
Also note that $H'(x,t)$ is not unimodal in $t$, but fixing $x$ and
restricting $t$ to only even or only odd values it becomes unimodal. As
$s_i(y)\sim y$ we can still apply Lemma~\ref{lunimod1} below.

Using (\ref{eq:inf}) and Lemmas~\ref{lunimod1} and~\ref{lhuni} we have
  \begin{align*}
    |\CON(y)| &= \bigg|\sum_{i \ge 0} (-1)^{i} \sum_{x \in X} \INF(y-x,t-s_i(y))\bigg|\\
    &= \bigg|\tfrac 12 \sum_{i \ge 0} (-1)^{i} \sum_{x \in X} \big[H(y-x+1,t-s_i(y)-1) \\
    &  \qquad\qquad - \; H(y-x-1, t-s_i(y)-1)\big]\bigg|\\
    &= \bigg|\tfrac 12 \sum_{i \ge 0} (-1)^{i} \big [H'(y+L,t-s_i(y)-1)\\
    &  \qquad\quad - \; H'(y, t-s_i(y)-1)\big]\bigg|\\
    &\le \bigg|\tfrac 12 \sum_{i \ge 0} (-1)^{i} H'(y+L,t-s_i(y)-1) \bigg| \\
    &  \qquad\quad + \; \bigg|\tfrac 12 \sum_{i \ge 0} (-1)^{i} H'(y, t-s_i(y)-1)\bigg|\\
    &\le \tfrac 12 \max_{s \in \N}   H'(y+L,s) + \tfrac 12 \max_{s \in \N}  H'(y, s)\\
    &=  O(1/(y+L-1/2))+O(1/(y-1/2)).
  \end{align*}

  Thus the vertices in $[-2L+1..L]$ contribute at most
  $$
  \sum_{y \in [-2L+1..L]} |\CON(y)| = O(\sum_{i = 1}^{2L}1/(i-1/2))= O(\log L).
  $$ 
  
  Combining all cases, we have \[|f(X,t) - E(X,t)| \le \sum_{y \in \Z}
  |\CON(y)| = O(\log L).\]
  
  For the lower bound, we just have to place the chips in a way
  the logarithmic contribution actually occurs. Without loss of
  generality, let $L$ be odd.

  Consider the following initial configuration (its existence is
  ensured by the parity forcing theorem): All arrows point
  towards the interval $X$ (arrows of vertices in $X$ may point
  anywhere). Let $t =  L^2$. Choose an initial configuration of the
  chips such that $f(y,s)$ is odd if and only if $y \in [L]$ is even
  and $t-s = y^2$.
  
  Now by construction, $\CON(y) = 0$ for all $y \in \Z \setminus [L]$.
  For $y \in [L]$, we have
  \begin{align*}
    \CON(y) &= - \sum_{x \in X} \INF(y-x,y^2)\\
    &= \tfrac 12 H(y,y^2) - \tfrac 12 H(y+L+1,y^2)\\
    &\ge \tfrac 12 H(y,y^2) - \tfrac 12 H(y+L+1,(y+L+1)^2)\\
    &= \Omega(y^{-1}).
  \end{align*}
  Hence for this initial configuration, \[f(X,t) - E(X,t) = \sum_{y
    \in \Z} \CON(y) = \sum_{y \in [L], y \sim 2} O(y^{-1}) = \Omega(\log L).\]
\end{proof}

\section{Intervals in Time}\label{sectint}

In this section, we regard the discrepancy in time-intervals.  For $x
\in \Z$ and finite $S \subseteq \N_0$, set
\begin{eqnarray*}
  f(x,S) &:=& \sum_{t \in S} f(x,t),\\
  E(x,S) &:=& \sum_{t \in S} E(x,t).
\end{eqnarray*}

We show that the discrepancy of a single vertex in a time-interval of
length $T$ is $O(\sqrt T)$, and this is sharp. 

\begin{theorem}\label{ttime}
  The maximal discrepancy $|f(x,S) - E(x,S)|$ of a single vertex~$x$
  in a time interval $S$ of length $T$ is $\Theta(T^{1/2})$.
\end{theorem}

In the proof, we need the following fact that ``rolling sums'' of
unimodal functions are unimodal again.

\begin{lemma}[Unimodality of rolling sums]\label{lunimod2}
  Let $f : \Z \to \R$ be unimodal. Let $k \in \N$. Define $F : \Z \to \R$ by
  $F(z) = \sum_{i=0}^{k-1} f(z+i)$. Then $F$ is unimodal.
\end{lemma}

\begin{proof}
  Let $f$ and $m \in \Z$ be such that $f$ is
  non-decreasing in $\Z_{\le m}$ and non-increasing in $\Z_{\ge m}$. We show
  that for some $m-k<M\le m$ we have that $G(x) := F(x+1) - F(x)$ is
  nonnegative for $x<M$ and nonpositive for $x\ge M$. This implies that $F$
  is unimodal.
  
  Since $G(x) = f(x+k) - f(x)$ for all $x \in \Z$, $G(x)$ is non-negative
  for $x\le m-k$ and it is nonpositive for $x\ge m$. For $m-k\le x<m$ we
  have $G(x+1) - G(x) = (f(x+k+1) - f(x+k)) - (f(x+1) - f(x)) \le
  0$, that is, $G$ is non-increasing in $[m-k..m]$. Hence $M$ exists as
  claimed.
\end{proof}

Of course, analogous statements hold for functions defined only on
even or odd integers.

The following result says that a single odd split has an influence of
exactly one on another vertex over infinite time.

\begin{lemma}\label{linftime}
   For all $x \in \Z \setminus \{0\}$, $\sum_{t \in \N} |\INF(x,t)| = 1$.
\end{lemma}

\begin{proof}
  W.l.o.g., let $x \in \N$. Then $|\INF(x,t)| = \tfrac 12 H(x-1,t-1) -
  \tfrac 12 H(x+1,t-1)$. Consider a random walk of a single chip
  started at zero. Let $X_{y,t}$ be the indicator random variable for
  the event that the chip is on vertex $y$ at time $t$. Let $Y_{y,t}$
  be the indicator random variable for the event that the chip is on
  vertex $y$ at time $t$ and that it has not visited vertex $x$ so
  far. Let $T$ denote the first time the chip arrives at $x$.
   
   For any $t>s>0$ we have by symmetry that
   $\Pr(X_{x-1,t-1}=1|T=s)=\Pr(X_{x+1,t-1}=1|T=s)$. Clearly, for $t \le
   T$, $X_{x+1,t-1} = 0$, and for $t > T$, $Y_{x-1,t-1} = 0$. Thus
   \begin{eqnarray*}
     \sum_{t \in \N} |\INF(x,t)| &=& \tfrac 12 \sum_{t \in \N}
     (E(X_{x-1,t-1}) - E(X_{x+1,t-1}))\\
     &=& \tfrac 12 \sum_{s \in \N} \Pr(T = s) \sum_{t \in \N}
     E((X_{x-1,t-1} - X_{x+1,t-1}) \setmid T = s)\\
     &=& \tfrac 12 \sum_{s \in \N} \Pr(T = s) \sum_{t \in [s]}
     E(X_{x-1,t-1} \setmid T = s)\\
     &=& \tfrac 12 \sum_{s \in \N} \Pr(T = s)  
     E\bigg(\sum_{t \in [s]} X_{x-1,t-1} \setmid T = s\bigg)\\
     &=& \tfrac 12 \sum_{s \in \N} \Pr(T = s)  
     E\bigg(\sum_{t \in \N} Y_{x-1,t-1} \setmid T = s\bigg)\\
     &=& \tfrac 12 E\bigg(\sum_{t \in \N} Y_{x-1,t-1}\bigg).     
   \end{eqnarray*}
   Note that $E(\sum_{t \in \N} Y_{x-1,t-1})$ is just the expected
   number of visits to $x-1$ before visiting $x$. This number of
   visits is exactly $k$ if and only if the chip moves left after each
   of its first $k-1$ visits and right after the $k$th visit. This
   happens with probability $2^{-k}$. Hence $E(\sum_{t \in \N}
   Y_{x-1,t-1}) = \sum_{i \in \N} i 2^{-i} = 2$.
\end{proof}

\begin{proof}[Proof of Theorem~\ref{ttime}]
  Fix any even initial configuration. Let $t_0 \in \N_0$ and $S = [t_0
  \,..\,t_0+T-1]$. Without loss, let $x = 0$. By (\ref{eqn:disc}), we
  have
  \begin{eqnarray*}
    f(0,S) - E(0,S) &=& \sum_{t \in S} (f(0,t) - E(0,t)) \\
    &=& \sum_{y \in \Z} \ARR(y,0) \sum_{i \ge 0} (-1)^{i} \sum_{t \in S} 
    \INF(y,t-s_i(y)).
  \end{eqnarray*}
  By unimodality of rolling sums (Lemma~\ref{lunimod2}),
  \begin{eqnarray*}
    |f(0,S) - E(0,S)| &\le& \sum_{y \in \Z} \,\, \max_{s \in \N} \,\, \bigg|\sum_{t \in S}   
    \INF(y,t-s)\bigg|.    
  \end{eqnarray*}
  We estimate the term $\max_{s \in \N} \,\, |\sum_{t \in S}
  \INF(y,t-s)|$ for all $y$.  For $1 \le |y| \le T^{1/2}$, we use
  Lemma~\ref{linftime} and simply estimate
  \begin{equation}
    \max_{s \in \N} \bigg|\sum_{t \in S} \INF(y,t-s)\bigg| \le \sum_{t
    \in \N} |\INF(y,t)| = 1.
  \end{equation}
  For $|y| > T^{1/2}$,
  $$
  \max_{s \in \N} \bigg|\sum_{t \in S} \INF(y,t-s)\bigg| \le T \max_{t \in \N} |\INF(y,t)| = T O(y^{-2})
  $$
  by Lemma~\ref{linfmax}.  Hence \[|f(0,S) - E(0,S)| \le \sum_{1
    \le |y| \le T^{1/2}} 1 \; + T \!\!\! \sum_{|y| > T^{1/2}} O(y^{-2}) =
  O(T^{1/2}).\]

  For the lower bound, we invoke the parity forcing theorem again. By
  this, there is an  even initial configuration such that all arrows
  point towards zero, and such that there is an odd number of chips on vertex $x
  \in \Z$ at time $t \in \N_0$ if and only if $x \in X := [\sqrt T \,..\, 2
  \sqrt T]$ and $t = 4T - x^2$. For this initial configuration and
  $S = [4T+1\,..\,5T]$, we compute
  \begin{eqnarray*}
    \lefteqn{|f(0,S) - E(0,S)|}\\ &=& \sum_{t \in S} \sum_{y \in X} |\INF(y,t
     - 4T + y^2)| \\
     &\ge& (1/2) T^{3/2} \min\big\{|\INF(y,t)| \,\big|\, y \in
    X, t \in S, y \sim t\big\} \\
    &=& \Omega(T^{1/2}).
  \end{eqnarray*}
  \end{proof}

\section{Space-Time-Intervals}\label{secstint}\label{sectspace}

We now regard the discrepancy in space-time-intervals. Extending the
previous notation, for finite $X \subseteq \Z$ and finite $S \subseteq
\N_0$ set 
\begin{eqnarray*}
  f(X,S) &:=& \sum_{x \in X} \sum_{t \in S} f(x,t),\\
  E(X,S) &:=& \sum_{x \in X} \sum_{t \in S} E(x,t).
\end{eqnarray*}

\begin{theorem}\label{tspacetime}
  Let $X \subseteq \Z$ and $S \subseteq \N_0$ be finite intervals of
  lengths $L$ and $T$, respectively. Then the maximal discrepancy
  $|f(X,S) - E(X,S)|$ (taken over all odd or even initial
  configurations) is $\Theta(T \log(L T^{-1/2}))$, if $L \ge 2 T^{1/2}$,
  and $\Theta(L T^{1/2})$ otherwise.
\end{theorem}

\begin{proof}
For the upper bound we use Theorems~\ref{tspace} and~\ref{ttime}.
To prove $|f(X,S) - E(X,S)|=O(LT^{1/2})$ we can simply apply
Theorem~\ref{ttime}:
  \begin{eqnarray*}
    |f(X,S) - E(X,S)| &\le& \sum_{x \in X} |f(x,S) - E(x,S)| \\
    &\le& L O(T^{1/2}).
  \end{eqnarray*}

For the other upper bound $|f(X,S) - E(X,S)|=O(T\log(LT^{-1/2}))$ we have to
separate contributions of the vertices and apply
the bounds in the proof of Theorem~\ref{tspace} for most of them and
the bounds from the proof of Theorem~\ref{ttime} for the rest.

  Fix an even initial configuration. Without loss of generality, let
  $X = [-L+1..0]$. Let $t_0 \in \N_0$ and $S = [t_0 \, .. \,
  t_0+T-1]$. As in previous proofs, by (\ref{eqn:disc}) we have $f(X,S)
  - E(X,S) = \sum_{y \in \Z} \CON(y)$ with \[\CON(y) := \ARR(y,0)
  \sum_{i \ge 0} (-1)^i \sum_{x \in X} \sum_{t \in S}
  \INF(y-x,t-s_i(y)).\]
  
Here $\CON(y)$ is the sum for $t\in S$ of the contribution $\CON_t(y)$ of $y$
to the discrepancy of the interval $X$ at a single time step $t$. The bound we
established in the proof of Theorem~\ref{tspace} is $|\CON_t(y)|=O(Ly^{-2})$
for $y>L$ and
$y\le-2L$ and $|\CON_t(y)|=O(1/(y-1/2)+1/(y+L-1/2))$ for $-2L<y\le L$. Thus
we have $$|\CON(y)|=O(LTy^{-2})$$ for $y>L$ and $y\le-2L$ and
$$|\CON(y)|=O(T/(y-1/2)+T/(y+L+1/2))$$ for $-2L<y\le L$.

The above bounds are the largest for $y$ close to $0$ or $-L$.
For $|y|\le T^{1/2}$ and for $|y+L|\le T^{1/2}$ we bound $|\CON(y)|$
in a different way. Let $X'$ be the interval
$[-L+2\lceil T^{1/2}\rceil..-2\lceil T^{1/2}\rceil]$ or empty if
$-L+2\lceil T^{1/2}\rceil>-2\lceil T^{1/2}\rceil$. We express the
contribution $\CON(y)$ of $y$ as the sum of contributions to
different parts of $X$. Let $\CON'(y)$ be
the total contribution of the vertex $y$ to the discrepancy in $X'$ over the
time interval $S$. Since $y$ is separated from $X'$ by at least $T^{1/2}$ the
above bound gives $\CON'(y)=O(T^{-1/2})$. Let $\CON'_x(y)$ be the total
contribution of $y$ to the discrepancy of the single vertex
$x\in X\setminus X'$ over the time interval $S$. To bound $\CON'_x(y)$
we apply the technique of the proof of Theorem~\ref{ttime}: by
Lemma~\ref{linftime} we have $|\CON'_x(y)|<1$. Thus we have
$$|\CON(y)|\le|\CON'_x(y)|+\sum_{x\in X\setminus X'}|\CON'_x(y)|=
O(T^{1/2})+O(T^{1/2})=O(T^{1/2}).$$

Let $H_1$ be the set of vertices $y$ with $y\le-2L$ or $y>L$. The total
contribution of these vertices is at most
$$\sum_{y\in H_1}|\CON(y)|=\sum_{y\in H_1}O(LTy^{-2})=O(T).$$
Let $H_2$ be the set of vertices $y$ with $|y|\le T^{1/2}$ or
$|y+L|\le T^{1/2}$. The total contribution of these vertices is at most
$$\sum_{y\in H_2}|\CON(y)|=\sum_{y\in H_2}O(T^{1/2})=O(T).$$
Let $H_3$ be the the set of vertices $y$ outside $H_1$ and $H_2$. Their
total contibution is bounded by
\begin{align*}
\sum_{y\in H_3}|\CON(y)| &=\sum_{y\in H_3}O(T/y+T/(y+L)) \\
&= O \left (T \!\!\!\! \sum_{i=\lceil T^{1/2}\rceil}^{2L}1/i \right)=O(T\log(L/T^{1/2})).
\end{align*}
Finally we have
\begin{align*}
|f(X,S)-E(X,S)|& =|\sum_{y \in \Z} \CON(y)| \\
& \le O(T)+O(T)+O(T\log(L/T^{1/2})) \\
& = O(T\log(L/T^{1/2})).
\end{align*}
  
  We now prove the corresponding lower bounds. Assume first that
$L\ge 2T^{1/2}$. Set $Y =
  [T^{1/2}..\,L]$. Choose an even initial configuration such that
  $f(x,t)$ is odd if and only if $x \in Y$ and $t = L^2 - x^2$.
  Direct all arrows towards zero. Let $S = [L^2 .. \,L^2 + T - 1]$.
  Then for $y \in Y$, with appropriately chosen $\delta_t, \ee_t \in \{0,1\}$
  we have
  \begin{eqnarray*}
    \CON(y) &=& \sum_{x \in X} \sum_{t = L^2}^{L^2 + T-1} \left |\INF(y-x,t-(L^2 - y^2)) \right|\\
    &\ge& \tfrac 12  \sum_{t = y^2}^{y^2+ T - 1} (H(y-1+\ee_t,t-1) - H(y+L-1-\delta_t,t-1))\\
    &\ge& \Omega\bigg(\sum_{t = y^2}^{y^2+ T - 1} H(y-1+\ee_t,t-1)\bigg)\\
    &=& \Omega(T y^{-1}).
  \end{eqnarray*}
  For $y \notin Y$, $\CON(y) = 0$. Hence the discrepancy in this
  setting is
  $$
  \sum_{y \in Y} \CON(y) = \sum_{y = T^{1/2}}^L \Omega(T
  y^{-1}) = \Omega(T \log(LT^{-1/2})).
  $$

Assume now that $L \le 2 T^{1/2}$. The setting of Theorem~\ref{ttime} works
for this lower bound, too. Choose an initial configuration such that $f(x,t)$
is odd if and only if $x \in X := [T^{1/2} .. \, 2
  T^{1/2}]$ and $t = 4T - y^2$. Then 
  \begin{align*}
    \CON(y) &= \sum_{x \in X} \sum_{t = 4T}^{5T-1} |\INF(y-x,t-(4T -
    y^2))| \\
    &\ge L T \min \left \{|\INF(y,t)| \; \bigg | \; y\in \Z \cap
    [T^{1/2} .. \, 3 T^{1/2}], t \in \Z \cap [T .. 5T], y \sim t \right \} \\
    &= \Omega(L)
  \end{align*}
  for all $y \in Y$. Again, $\CON(y) = 0$ for $y \notin Y$.  Hence
  $\sum_{y \in \Z} \CON(y) = \Omega(L T^{1/2})$.
\end{proof}

\section{Intervals in Space, Revisited}

We stated in Theorem~\ref{intspace} that the discrepancy in an interval of
length $L$ is $O(\log L)$. Here we show that intervals of length $L$ with
about $\log L$ discrepancy are very rare, the root-mean-squared (i.e., quadratic) average of the
discrepancies of a
long contiguous set of intervals of length $L$ is only $O(\sqrt{\log L})$, and
this bound is tight.
 
For a set $X$ of vertices we denote by $\DISC(X,t)$ the discrepancy of the set
$X$ at time $t$, i.e., we set $\DISC(X,t)=f(X,t)-E(X,t)$.

\begin{theorem}
  Let $X$ be an interval of length $L$. For $M$ sufficiently large,
  \[\frac 1 M \sum_{k = 1}^M\DISC^2(X+k,t) = O(\log L).\]
Furthermore, for a given $L$ and $M$ there exists an even initial
configuration, and a time $t$ and an interval $X$ of length $L$ such that
$$\frac1M\sum_{k=1}^M\DISC^2(X+k,t)=\Omega(\log L).$$
\end{theorem}

\begin{proof}
For the first statement we need to prove an $O(\sqrt{\log L})$ bound on the
quadratic average of the discrepancies $\DISC(X+k,t)$ with $k=1,\ldots,M$.
First note that by changing the
individual discrepancies by a bounded amount, we change the
quadratic average by at most the same amount. We use this observation to freely
neglect $O(1)$ terms in the discrepancy of the intervals. In particular we can
change the intervals themselves by adding or deleting a
bounded number of vertices. We use this to make a few simplifying assumptions.
As in Section~\ref{sectspace} we assume that (i) the starting configuration is
odd, (ii) the interval $X$ is $X=[-L'..L']$ with $L'\sim t$, and (iii) $M$ is
even and we only consider even values of $k$, i.e., we consider the average of
$\DISC^2(X+k,t)$ for $2\le k\le M$, $k$ even (this can be justified
by considering $X+k+1$ instead of $X+k$ for odd $k$).

First we show that discrepancies caused by odd piles at time $t-L^2$ or before
can be neglected. We start with (\ref{eqn:disc}) for the
individual discrepancies $\DISC(x,t)$.
\begin{eqnarray*}
\DISC(x,t)&=&\sum_{y\in\Z}\ARR(y,0)\sum_{i\ge0}(-1)^i\INF(y-x,t-s_i(y))\\
&=&\DISC_1(x,t)+\DISC_2(x,t);\\
\DISC_1(x,t)&=&
\sum_{y\in\Z}\ARR(y,0)\sum_{s_i(y)>t-L^2}(-1)^i\INF(y-x,t-s_i(y));\\ 
|\DISC_2(x,t)|&=&
\left | \sum_{y\in\Z}\ARR(y,0)\sum_{s_i(y)\le t-L^2}(-1)^i\INF(y-x,t-s_i(y)) \right |\\
&\le&\sum_{y\in\Z}\max_{u\ge L^2}|\INF(y-x,u)|.
\end{eqnarray*}
We have seen that  $|\INF(z,u)|$ is unimodal for fixed $z$ and its maximum is
at $u=\lceil z^2/3\rceil$, so we have
\begin{eqnarray*}
|\DISC_2(x,t)|&\le&
2\sum_{z=1}^L|\INF(z,L^2)|+2\sum_{z>L}|\INF(z,\lceil z^2/3\rceil)|\\
&\le&2\sum_{z=1}^L\frac z{L^2}H(z,L^2)+2\sum_{z>L}O(z^{-2})\\
&\le&2\sum_{z=1}^LH(z,L^2)/L+O(1/L)=O(1/L).
\end{eqnarray*}
Therefore the total contribution of $\DISC_2$ to the discrepancy of an
interval $X+k$ is small. For
$$\DISC_1(X+k,t):=\sum_{x\in X+k}\DISC_1(x,t)$$
we have
$$|\DISC_1(X+k,t)-\DISC(X+k,t)|=\left | \sum_{x\in X+k}\DISC_2(x,t) \right |=O(1).$$

We continue as in Section~\ref{sectspace} collapsing a sum using
$\INF(z,u)=\frac12H(z+1,u-1)-\frac12H(z-1,u-1)$. We also use that
$\DISC(x,t)=\DISC_1(x,t)=0$ for $x\sim t$ as the starting configuration is
odd.
\begin{eqnarray*}
&&\DISC_1(X+k,t)\\&=&\sum_{x\in X+k,x\sim t+1}\sum_{y\in\Z}
\ARR(y,0)\sum_{s_i(y)>t-L^2}(-1)^i\INF(y-x,t-s_i(y))\\
&=&\sum_{y\in\Z}\ARR(y,0)\sum_{s_i(y)>t-L^2}(-1)^i\sum_x
\left(\frac12H(y-x+1,t-s_i(y)-1)\right.\\&&\left.-\frac12H(y-x-1,t-s_i(y)-1)\right)\\
&=&\frac12\sum_{y\in\Z}\ARR(y,0)\sum_{s_i(y)>t-L^2}(-1)^i
\left(H(y-k+L',t-s_i(y)-1)\right.\\&&\left.-H(y-k-L',t-s_i(y)-1)\right).
\end{eqnarray*}
We separate the two terms in this last expression. With
$$D(m):=2\sum_{y\in\Z}\ARR(y,0)\sum_{s_i(y)>t-L^2}(-1)^iH(y-m,t-s_i(y)-1)$$
we have
$$\DISC_1(X+k,t)=\frac14D(k-L')-\frac14D(k+L').$$
Our original goal was to prove an $O(\sqrt{\log L})$ bound on the quadratic
average of $\DISC(X+k,t)$. As $\DISC_1(X+k,t)$ differs from $\DISC(X+k,t)$ by
$O(1)$ it is clearly enough to prove the same bound for the quadratic average
of $\DISC_1(X+k,t)$. By the last displayed formula it is enough to prove the
$O(\sqrt{\log L})$ bound on the two parts
$D(k-L')$ and $D(k+L')$ separately, both for $0<k\le M$ even. It is therefore
enough to bound the quadratic average of $D(m)$ for an arbitrary interval $I$
of length $M$. Here we consider only values $m\sim t$, for other values of $m$
we have $D(m)=0$.

Let $t_0=\max(0,t-L^2+1)$ be the first time-step considered. For $y\in\Z$ and
$u\sim y+1$ we have an odd pile at $y$ if and only if
$\ARR(y,u)\ne\ARR(y,u+2)$ and in this case $\ARR(y,u)=(-1)^i\ARR(y,0)$ for the
index $i$ with $s_i(y)=u$. We estimate the contribution $D(m,y)$ of a fixed
value $y$ to the sum defining $D(m)$. For $m\sim t$ we have
\begin{eqnarray*}
&&D(m,y):= 2\ARR(y,0)\sum_{s_i(y)>t-L^2}(-1)^iH(y-m,t-s_i(y)-1)\\&=&
\sum_{t_0\le u<t, u\sim y+1}(\ARR(y,u)-\ARR(y,u+2))H(y-m,t-u-1)\\
&=&\sum_{t_0+2\le u<t, u\sim y+1}\ARR(y,u)(H(y-m,t-u-1)-H(y-m,t-u+1))\\&&+
\ARR(y,t_1)H(y-m,t-t_1-1)-\ARR(y,t_2)H(y-m,t-t_2+1),
\end{eqnarray*}
where $t_1=t_1(y)$ is either $t_0$ or $t_0+1$, whichever makes $t_1\sim y+1$
and similarly $t_2=t_2(y)$ is either $t$ or $t+1$, so that $t_2\sim
y+1$. We have
$$D(m)=\sum_{y\in\Z}D(m,y)$$
and with
$$
D'(m):= \sum_{y\in\Z} \; \sum_{\substack{t_0+2\le u<t-2 \\ u\sim y+1}} \ARR(y,u)(H(y-m,t-u-1)-H(y-m,t-u+1))
$$
we have
\begin{align*}
|D(m)-D'(m)| &= \bigg |\sum_{y\in\Z} (\ARR(y,t_1(y))H(y-m,t-t_1(y)-1) \\
& \qquad \quad - \;\ARR(y,t_2(y))H(y-m,t-t_2(y)+1)) \bigg |\\
&\le \sum_{u\in\{0,1,t-t_0-2,t-t_0-1\}} \; \sum_{y\in\Z}H(y-m,u)\le4.
\end{align*}

As before, we ignore the small difference and will prove the $O(\sqrt{\log
L})$ bound on the quadratic average of $D'(m)$ instead of $D(m)$. Computing
the square and summing over $m$ we get the following. The
summations are taken for $m\in I$, $m\sim t$, for $y_,y_2\in\Z$, and for
$u_1,u_2\in[t_0+2,t-3]$, $u_1\sim u_2\sim y+1$, respectively.
\begin{eqnarray*}
\sum_mD'^2(m)&=&\sum_{y_1,y_2}\sum_{u_1,u_2}
\ARR(y_1,u_1)\ARR(y_2,u_2)\\&&\cdot\sum_m
(H(y_1-m,t-u_1-1)-H(y_1-m,t-u_1+1))\\&&\cdot
(H(y_2-m,t-u_2-1)-H(y_2-m,t-u_2+1)).
\end{eqnarray*}

Let us estimate the contribution to this sum coming from a fixed $y_1$, $u_1$,
and $u_2$. Disregarding the signs and extending the summation for all $m$ (even
outside $I$) the contribution of each of the four terms we get from the
multiplication is exactly $1$. As $u_1$ and $u_2$ can take at most $L^2/2$
values each, the total contribution coming from a single value of $y_1$ is at
most $L^4$.

Let us obtain the intervals $I'$ and $I''$ from $I$ by extending or shortening
it at both ends by $L^2$ respectively, i.e., if $I=[a,b]$, then $I'=[a-L^2,b+L^2]$,
$I''=[a+L^2,b-L^2]$. If $y_1$ is outside $I'$ we have
$H(y_1-m,t-u_1-1)=H(y_1-m,t-u_1+1)=0$ for all $m\in I$, therefore such $y_1$
has zero contribution to $\sum_mD'^2(m)$. The contribution for fixed $y_1$,
$y_2$, $u_1$, and $u_2$ can usually be written in closed form
using the identity
$$\sum_mH(y_1-m,v_1)H(y_2-m,v_2)=H(y_1-y_2,v_1+v_2).$$
This identity is valid if we sum over all possible values of $m$, but for $y_1\in
I''$ the contribution of the values $m\notin I$ is zero. Therefore the
contribution to $\sum_mD'^2(m)$ of the fixed terms $y_1\in I''$, $y_2$,
$u_1$, and $u_2$ is
\begin{align*}
& \ARR(y_1,u_1) \ARR(y_2,u_2) \sum_m (H(y_1-m,t-u_1-1)-H(y_1-m,t-u_1+1)) \\
& \qquad \quad \cdot (H(y_2-m,t-u_2-1)-H(y_2-m,t-u_2+1)) \\
&= \ARR(y_1,u_1)\ARR(y_2,u_2)(H(y,v-2)-2H(y,v)+H(y,v+2)),
\end{align*}
where $y=y_1-y_2$ and $v=2t-u_1-u_2$.

To estimate these contributions we first calculate
$$H(y,v-2)-2H(y,v)+H(y,v+2)=O(y^4/v^4+1/v^2)H(y,v+2).$$
The same $y=y_1-y_2$ value arises exactly once for every $y_1\in I''$, a
total of $M-2L^2$ possibilities. The largest possible value of $v$ is less than
$2L^2$ and any single value $v$ can be the result of at most $v$ pairs $u_1$,
$u_2$. There are $4L^2$ possible values of $y_1$ outside $I''$ but inside $I'$
contributing at most $4L^6$. Summing for all these contributions we estimate
\begin{eqnarray*}
\sum_mD'^2(m)&\le&4L^6+O\left(\sum_{v=1}^{2L^2}Mv\sum_{y\in\Z}
(y^4/v^4+1/v^2)H(y,v+2)\right)\\
&=&4L^6+O\left(M\sum_{v=1}^{2L^2}\sum_{y\in\Z}(y^4/v^3+1/v)H(y,v+2)\right)\\
&=&4L^6+O\left(M\sum_{v=1}^{2L^2}1/v\right)=O(L^6+M\log L).
\end{eqnarray*}
Here we used the estimate on the fourth moment of the random walk:
$$\sum_{y\in\Z}y^4H(y,v+2)=O((v+2)^2)=O(v^2).$$

To finish the proof we set the threshold $M>L^6$ for sufficiently large $M$.
We did not make an effort to optimize for this threshold.
This ensures that $\sum_mD'^2(m)=O(M\log L)$, so the quadratic average of
$D'(m)$ (and therefore of $\DISC(X+k,t)$) is $O(\sqrt{\log L})$ as claimed.

It remains to construct a starting configuration where the quadratic average of
discrepancies in the intervals of length $L$ is large. For our construction we
do not even use the value $L$. For a given (even) parameter $t$, we define a
probability distribution on starting positions, such that for all $L<t$
and all intervals $X$ of length $L$ the expectation of
$\DISC^2(X,t)=\Omega(\log L)$.

We let $r(a,b)$ stand for independent random $\pm1$ variables for all integers
$a$ and $b\ge1$. We look for an even starting configuration (guaranteed by the
Arrow-Forcing Theorem), such that $\ARR(x,u)=r(a,b)$ for all even $x$ and
$u$ satisfying $4^b<u\le4^{b+1}$ and $a2^b<x\le(a+1)2^b$. 
For simplicity we set $\ARR(x,u)=1$ for all $u$ and all odd $x$ and we
also set $\ARR(x,u)=1$ for all $x$ and $u\le4$.

A simple calculation similar to the one in Section~\ref{sectspace} shows that
for an interval $X=[c,d]$ we have
$$\DISC(X,t)=\sum_{a,b}h(a,b)r(a,b),$$
where the coefficients $h(a,b)$ depend on $X$. Further analysis shows that all
coefficients are bounded and $\Theta(\log L)$ of them are above a positive
absolute constant for each interval of length $L$. This implies that the
expectation of $\DISC^2(X,t)$ is $\Omega(\log L)$, and therefore the
expectation of the average $\frac1M\sum_{k=1}^M\DISC^2(X+k,t)$ is also
$\Omega(\log L)$. This proves the second statement of the theorem.
\end{proof}


\begin{thebibliography}{Kle05}
\label{sec:biblio}

\bibitem[CS04]{cooperproppx}
J.~Cooper and J.~Spencer.
\newblock {Simulating a Random Walk with Constant Error}.
\newblock {\tt arXiv:math.CO/0402323}.

\bibitem[CS05]{cooperpropp}
J.~Cooper and J.~Spencer.
\newblock Simulating a random walk with constant error.
\newblock {\em Combinatorics, Probability and Computing}.
\newblock To appear.

\bibitem[Kle05]{kleber}
M.~Kleber.
\newblock Goldbug variations.
\newblock {\em Mathematical Intelligencer}, 27:55--63, 2005.

\bibitem[LP05]{levineperes}
L.~Levine and Y.~Peres.
\newblock {Spherical Asymptotics for the Rotor-Router Model in $\mathbb{Z}^d$}.
\newblock {\tt arXiv:math.PR/0503251}.

\end{thebibliography}
\end{document}